\documentclass[12pt,twoside]{amsart}

\usepackage{amsopn,amssymb}
\usepackage{tikz}
\usepackage[margin=1in]{geometry}
\usepackage[pdfauthor={David Cook II, Sonja Mapes, Gwyneth Whieldon},
            pdftitle={Partially ordered sets in Macaulay2},
            pdfproducer={LaTeX with hyperref},
            pdfcreator={latex->dvips->ps2pdf},
            pdfborder={0 0 0},
            colorlinks=true]{hyperref}

\newcommand{\rank}{\text{rank}}
\newcommand{\cA}{\mathcal{A}}
\newcommand{\cH}{\mathcal{H}}
\newcommand{\cM}{\mathcal{M}}
\newcommand{\tH}{\text{H}}
\newcommand{\CC}{\mathbb{C}}
\newcommand{\RR}{\mathbb{R}}
\newcommand{\ZZ}{\mathbb{Z}}

\renewcommand{\section}[2] {\vspace{0.75\baselineskip} \noindent {\scshape #2}.}

\newtheorem*{theorem*}{Theorem}
\theoremstyle{definition}

\begin{document}

\title[Posets in \emph{Macaulay2}]{Partially ordered sets in \emph{Macaulay2}}

\author[D.\ Cook II]{David Cook II}
\address{Department of Mathematics, University of Notre Dame, Notre Dame, IN 46556, USA}
\email{\href{mailto:dcook8@nd.edu}{dcook8@nd.edu}}

\author[S.\ Mapes]{Sonja Mapes}
\email{\href{mailto:smapes1@nd.edu}{smapes1@nd.edu}}

\author[G.\ Whieldon]{Gwyneth Whieldon}
\address{Department of Mathematics, Hood College, Frederick, MD 21701, USA}
\email{\href{mailto:whieldon@hood.edu}{whieldon@hood.edu}}

\thanks{\emph{Posets} version 1.0.6 available at \url{http://www.nd.edu/~dcook8/files/Posets.m2}.}

\begin{abstract}
    We introduce the package \emph{Posets} for \emph{Macaulay2}.  This package provides a data structure
    and the necessary methods for working with partially ordered sets, also called posets.  In
    particular, the package implements methods to enumerate many commonly studied classes of posets,
    perform operations on posets, and calculate various invariants associated to posets.
\end{abstract}

\maketitle

\section*{Introduction}

A \emph{partial order} is a binary relation $\preceq$ over a set $P$ that is antisymmetric, reflexive, and transitive.
A set $P$ together with a partial order $\preceq$ is called a \emph{poset}, or \emph{partially ordered set}.

Posets are combinatorial structures that are used in modern mathematical research, particularly in algebra.
We introduce the package \emph{Posets} for \emph{Macaulay2} via three distinct posets or related ideals which
arise naturally in combinatorial algebra.  

We first describe two posets that are generated from algebraic objects.  The intersection semilattice associated
to a hyperplane arrangement can be used to compute the number of unbounded and bounded real regions cut out by
a hyperplane arrangement, as well as the dimensions of the homologies of the complex complement of a hyperplane arrangement.

Given a monomial ideal, the lcm-lattice of its minimal generators gives information on the structure of the free
resolution of the original ideal.  Specifically, two monomial ideals with isomorphic lcm-lattices have the
``same'' (up to relabeling) minimal free resolution, and the lcm-lattice can be used to compute, among other things,
the multigraded Betti numbers $\beta_{i,\bf{b}}(R/M)=\dim_{\Bbbk}\text{Tor}_{i,\bf{b}}(R/M,\Bbbk)$ of the monomial ideal.

In contrast to the first two examples (associating a poset to an algebraic object), we then describe an ideal that is
generated from a poset.  In particular, the Hibi ideal of a finite poset is a squarefree monomial ideal which has many
nice \emph{algebraic} properties that can be described in terms of \emph{combinatorial} properties of the poset.  In
particular, the resolution and Betti numbers, the multiplicity, the projective dimension, and the Alexander dual are
all nicely described in terms of data about the poset itself.

\section*{Intersection (semi)lattices}

A \emph{hyperplane arrangement} $\cA$ is a finite collection of affine hyperplanes in some vector space $V$.
The \emph{dimension} of a hyperplane arrangement is defined by $\dim(\cA)=\dim(V)$, and the \emph{rank} of a
hyperplane arrangement $\rank(\cA)$ is the dimension of the span in $V$ of the set of normals to the
hyperplanes in $A$.

The \emph{intersection semilattice} ${\mathcal L(A)}$ of $\cA$ is the set of the nonempty intersections
of subsets of hyperplanes $\bigcap_{\cH\in \cA'}\cH$ for $\cH\in\cA'\subseteq \cA$, ordered by reverse
inclusion.  We include the empty intersection corresponding to $\cA'=\emptyset$, which is the minimal
element in the intersection meet semilattice $\hat{0} \in {\mathcal L(A)}$.  If the intersection of all
hyperplanes in $\cA$ is nonempty, $\bigcap_{\cH\in \cA} H \neq 0$, then the intersection meet semilattice
${\mathcal L(A)}$ is actually a lattice.  Arrangements with this property are called \emph{central arrangements}.

\begin{figure}[!ht]
    \begin{minipage}[b]{0.48\linewidth}
        \centering
        \begin{tikzpicture}[scale=1,auto=left,vertices/.style={circle, fill=black, inner sep=1.5pt}]
            \draw[-] (-1.5,1.5)--(2,-2) {};
            \draw[-] (0,1.5)--(0,-2) {};
            \draw[-] (-2,-2)--(1.5,1.5) {};
            \draw[-] (-2,-1)--(2,-1) {};
            \node[vertices, label=right:{$p_1$}] at (0,0) {};
            \node[vertices, label=above:{$p_2$}] at (-1,-1) {};
            \node[vertices, label=below left:{$p_3$}] at (0,-1) {};
            \node[vertices, label=below:{$p_4$}] at (1,-1) {};
            \node at (-1.8,1.8) {$\cH_1$};
            \node at (0,1.8) {$\cH_2$};
            \node at (1.8,1.8) {$\cH_3$};
            \node at (2.3,-1) {$\cH_4$};
        \end{tikzpicture}\\
        \emph{(i) $\cA$}
    \end{minipage}
    \begin{minipage}[b]{0.48\linewidth}
        \centering
        \begin{tikzpicture}[scale=1, vertices/.style={draw, fill=black, circle, inner sep=0pt}]
             \node [vertices, label=below:{$\hat{0}$}] (8) at (-0+0,0){};
             \node [vertices, label=left:{$\cH_1$}] (0) at (-2.25+0,1.5){};
             \node [vertices, label=left:{$\cH_2$}] (1) at (-2.25+1.5,1.5){};
             \node [vertices, label=right:{$\cH_3$}] (3) at (-2.25+3,1.5){};
             \node [vertices, label=right:{$\cH_4$}] (4) at (-2.25+4.5,1.5){};
             \node [vertices, label=left:{$p_1$}] (2) at (-2.25+0,3.0){};
             \node [vertices, label=left:{$p_4$}] (5) at (-2.25+1.5,3.0){};
             \node [vertices, label=right:{$p_3$}] (6) at (-2.25+3,3.0){};
             \node [vertices, label=right:{$p_2$}] (7) at (-2.25+4.5,3.0){};
             \foreach \to/\from in {0/5, 0/2, 1/2, 1/6, 3/2, 3/7, 4/5, 4/6, 4/7, 8/0, 8/1, 8/3, 8/4}\draw [-] (\to)--(\from);
        \end{tikzpicture}\\
        \emph{(ii) ${\mathcal L(A)}$}
    \end{minipage}
    \caption{The non-central hyperplane arrangement 
    	\[
    		\cA=\{\cH_1=V(x+y),\cH_2=V(x),\cH_3=V(x-y),\cH_4=V(y+1)\}
    	\]
    	and its intersection semilattice ${\mathcal L(A)}$}
    \label{fig:hyp-arr}
\end{figure}
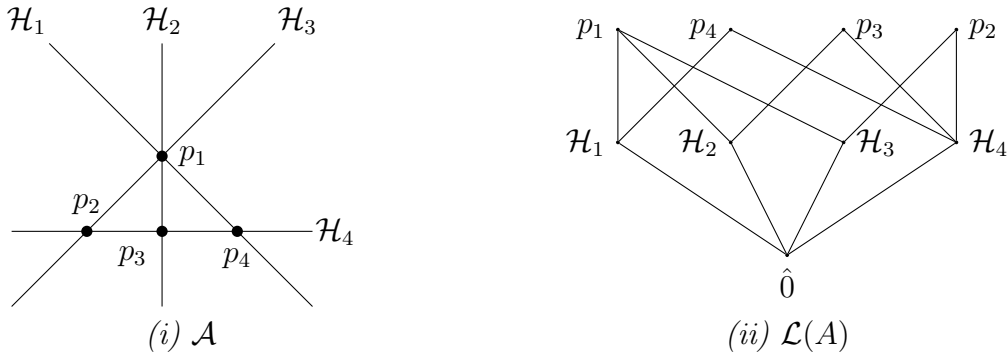

Consider the non-central hyperplane arrangement $\cA=\{\cH_1=V(x+y),\cH_2=V(x),\cH_3=V(x-y),\cH_4=V(y+1)\}$,
where $\cH_i=V(\ell_i(x,y))\subseteq \RR^2$ denotes the hyperplane $H_i$ of zeros of the linear form $\ell_i(x,y)$;
see Figure~\ref{fig:hyp-arr}(i).  We can construct ${\mathcal L(A)}$ in \emph{Macaulay2} as follows.  
\begin{verbatim}
i1 : needsPackage "Posets";
i2 : R = RR[x,y];
i3 : A = {x + y, x, x - y, y + 1};
i4 : LA = intersectionLattice(A, R);
\end{verbatim}
Further, using the method \texttt{texPoset} we can generate \LaTeX~to display the Hasse diagram of
${\mathcal L(A)}$, as in Figure~\ref{fig:hyp-arr}(ii).

A theorem of Zaslavsky~\cite{Zaslavsky} provides information about the topology of the complement of
hyperplane arrangements in $\RR^n$.  Let $\mu$ denote the M\"obius function of the intersection semilattice ${\mathcal L(A)}$.  Then
the number of regions that $\cA$ divides $\RR^n$ into is
\[
    r(\cA) = \sum_{x \in {\mathcal L(A)}} |\mu(\hat{0}, x)|.
\]
Moreover, the number of these regions that are bounded is
\[
    b(\cA) = |\mu({\mathcal L(A)} \cup \hat{1})|,
\]
where ${\mathcal L(A)} \cup \hat{1}$ is the intersection semilattice adjoined with a maximal element.

We verify these results for the non-central hyperplane arrangement $\cA$ using \emph{Macaulay2}:
\begin{verbatim}
i5 : realRegions(A, R)
o5 = 10
i6 : boundedRegions(A, R)
o6 = 2
\end{verbatim}

Moreover, in the case of hyperplane arrangements in $\CC^n$, using a theorem of Orlik and Solomon~\cite{OrlikSolomon},
we can recover the Betti numbers (dimensions of homologies) of the complement $\cM_\cA = \CC^n-\cup \cA$ of the hyperplane
arrangement using purely combinatorial data of the intersection semilattice.  In particular, $\cM_\cA$ has torsion-free integral
cohomology with Betti numbers given by
\[
    \beta_i({\cM}_\cA)=\dim_{\CC}\biggl(\tH_i({\cM}_\cA)\biggr)=\sum_{\substack{
        x\in{\mathcal L(A)} \\
        \dim^{\mathbb C}(x)=n-i
    }} |\mu(\hat{0},x)|,
\]
where $\mu(\cdot)$ again represents the M\"obius function.  See \cite{Wachs} for details and generalizations of this formula.

\emph{Posets} will compute the ranks of elements in a poset, where the ranks in the intersection lattice \texttt{LA} are
determined by the codimension of elements.  Combining the outputs of our rank function with the M\"obius function allows
us to calculate $\beta_0({\cM}_\cA) = 1$, $\beta_1({\cM}_{\cA}) = 4$, and $\beta_2({\cM}_{\cA}) = 5$.

\begin{verbatim}
i7 : RLA = rank LA
o7 =  {{ideal 0}, {ideal(x+y), ideal(x), ideal(x-y), ideal(y+1)}, 
       {ideal(y,x), ideal(y+1,x-1),ideal(y+1,x), ideal(y+1,x+1)}}
i8 : MF = moebiusFunction LA;
i9 : apply(RLA, r -> sum(r, x -> abs MF#(ideal 0_R, x)))
o9 = {1, 4, 5}
\end{verbatim}

\section*{LCM-lattices}

Let $R = K[x_1, \dots, x_t]$ be the polynomial ring in $t$ variables over the field $K$, where the degree of $x_i$ 
is the standard basis vector $e_i \in \ZZ^t$.  Let $M = (m_1, \dots, m_n)$ be a monomial ideal in $R$, then we define the
\emph{lcm-lattice} of $M$, denoted $L_M$, as the set of all least common multiples of subsets of the generators of
$M$ partially ordered by divisibility.  It is easy to see that $L_M$ will always be a finite atomic lattice.
While lcm-lattices are nicely structured, they can be difficult to compute by hand especially for large examples 
or for ideals where $L_M$ is not ranked.

\begin{figure}[!ht]
	\begin{tikzpicture}[scale=1, vertices/.style={draw, fill=black, circle, inner sep=0pt}]
		\node [vertices, label=right:{$1$}] (0) at (-0+0,0){};
		\node [vertices, label=right:{$a^{3} b^{2} c$}] (8) at (-3+0,1.33333){};
		\node [vertices, label=right:{$b^{2} c^{2} d$}] (1) at (-3+1.5,1.33333){};
		\node [vertices, label=right:{$a b c^{2} d$}] (2) at (-3+3,1.33333){};
		\node [vertices, label=right:{$a^{2} c d$}] (4) at (-3+4.5,1.33333){};
		\node [vertices, label=right:{$a^{3} b^{2} d$}] (7) at (-3+6,1.33333){};
		\node [vertices, label=right:{$a b^{2} c^{2} d$}] (3) at (-1.5+0,2.66667){};
		\node [vertices, label=right:{$a^{2} b c^{2} d$}] (5) at (-1.5+1.5,2.66667){};
		\node [vertices, label=right:{$a^{3} b^{2} c d$}] (9) at (-1.5+3,2.66667){};
		\node [vertices, label=right:{$a^{2} b^{2} c^{2} d$}] (6) at (-0+0,4){};
		\node [vertices, label=right:{$a^{3} b^{2} c^{2} d$}] (10) at (-0+0,5.33333){};
		\foreach \to/\from in {0/8, 0/1, 0/2, 0/4, 0/7, 1/3, 2/5, 2/3, 3/6, 4/5, 4/9, 5/6, 6/10, 7/9, 8/9, 9/10} \draw [-] (\to)--(\from);
	\end{tikzpicture}
    \caption{The lcm-lattice for $M = (a^3b^2c, a^3b^2d, a^2cd, abc^2d, b^2c^2d)$}
    \label{fig:lcm-lattice}
\end{figure}
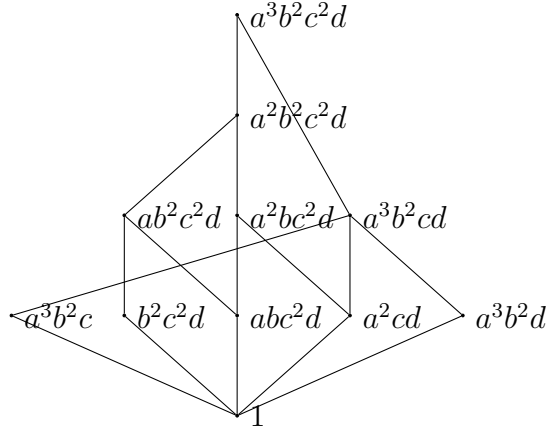

Consider the ideal $M = (a^3b^2c, a^3b^2d, a^2cd, abc^2d, b^2c^2d)$ in $R = k[a,b,c,d]$.  Then we can construct
$L_M$ in \emph{Macaulay2} as follows.  See Figure~\ref{fig:lcm-lattice} for the Hasse diagram of $L_M$, as
generated by the \texttt{texPoset} method.
\begin{verbatim}
i10 : R = QQ[a,b,c,d];
i11 : M = ideal(a^3*b^2*c, a^3*b^2*d, a^2*c*d, a*b*c^2*d, b^2*c^2*d);
i12 : LM = lcmLattice M;
\end{verbatim}

Lcm-lattices, which were introduced by Gasharov, Peeva, and Welker \cite{GasharovPeevaWelker}, have become an
important tool used in studying free resolutions of monomial ideals.  There have been a number of results
that use the lcm-lattice to give constructive methods for finding free resolutions for monomial ideals, for some examples see ,
\cite{Clark}, \cite{PeevaVelasco}, and \cite{Velasco}.

In particular, Gasharov, Peeva, and Welker~\cite{GasharovPeevaWelker} provided a key connection between the lcm-lattice
of a monomial ideal $M$ of $R$ and its minimal free resolution, namely, one can compute the (multigraded) Betti numbers
of $R/M$ using the lcm-lattice.  Let $\Delta(P)$ denote the order complex of the poset $P$, then for $i \geq 1$ we have
\[
	\beta_{i,b}(R/M) = \dim \tilde{H}_{i-2}(\Delta (\hat{0},b); k),
\]
for all $b \in L_M$, and so
\[
	\beta_i(R/M) = \sum_{b \in L_M} \dim \tilde{H}_{i-2}(\Delta (\hat{0},b); k).
\]

These computations can all be done using \emph{Posets} together with the package \emph{SimplicialComplexes},
by S.\ Popescu, G.\ Smith, and M.\ Stillman.  In particular, we can show that $\beta_{i, a^2b^2c^2d} = 0$
for all $i$ with the following calculation.

\begin{verbatim}
i13 : D1 = orderComplex(openInterval(LM, 1_R, a^2*b^2*c^2*d));
i14 : prune HH(D1)
o14 = -1 : 0
       0 : 0
       1 : 0
o14 : GradedModule
\end{verbatim}

Similarly, we can show that $\beta_{1, a^3b^2cd} = 2$.
\begin{verbatim}
i15 : D2 = orderComplex(openInterval(L, 1_R, a^3*b^2*c*d));
i16 : prune HH(D2)
o16 = -1 : 0
             2
       0 : QQ
o16 : GradedModule
\end{verbatim}

\section*{Hibi ideals}

Let $P = \{p_1, \ldots, p_n\}$ be a finite poset with partial order $\preceq$, and let $K$ be a field.
The \emph{Hibi ideal}, introduced by Herzog and Hibi~\cite{HibiHerzog}, of $P$ over $K$ is the squarefree
ideal $H_P$ in $R = K[x_1, \ldots, x_n, y_1, \ldots, y_n]$ generated by the monomials
\[
    u_I := \prod_{p_i \in I} x_i \prod_{p_i \notin I} y_i,
\]
where $I$ is an order ideal of $P$, i.e., for every $i \in I$ and $p \in P$, if $p \preceq i$, then $p \in I$.
\emph{Nota bene:} The Hibi ideal is the ideal of the monomial generators of the Hibi ring, a toric ring first
described by Hibi~\cite{Hibi}.

\begin{verbatim}
i17 : P = divisorPoset 12;
i18 : HP = hibiIdeal P; 
i19 : HP_*
o19 = {x x x x x x , x x x x x y , x x x x y y , x x x x y y , x x x y y y ,
        0 1 2 3 4 5   0 1 2 3 4 5   0 1 2 4 3 5   0 1 2 3 4 5   0 1 3 2 4 5
       x x x y y y , x x y y y y , x x y y y y , x y y y y y , y y y y y y }
        0 1 2 3 4 5   0 2 1 3 4 5   0 1 2 3 4 5  0 1 2 3 4 5   0 1 2 3 4 5
\end{verbatim}

Herzog and Hibi~\cite{HibiHerzog} proved that every power of $H_P$ has a linear resolution, and the 
$i^{\rm th}$ Betti number $\beta_i(R/H_P)$ is the number of intervals of the distributive
lattice $\mathcal{L}(P)$ of $P$ isomorphic to the rank $i$ boolean lattice.  Using Exercise~3.47 in
Stanley's book~\cite{Stanley}, we can recover this by looking instead at the number of elements of
$\mathcal{L}(P)$ that cover exactly $i$ elements.

\begin{verbatim}
i20 : betti res HP
             0  1  2 3
o20 = total: 1 10 12 3
          0: 1  .  . .
          5: . 10 12 3
i21 : LP = distributiveLattice P;
i22 : cvrs = partition(last, coveringRelations LP);
i23 : iCvrs = tally apply(keys cvrs, i -> #cvrs#i);
i24 : gk = prepend(1, apply(sort keys iCvrs, k -> iCvrs#k))
o24 : {1, 6, 3}
i25 : apply(#gk, i -> sum(i..<#gk, j -> binomial(j, i) * gk_j))
o25 : {10, 12, 3}
\end{verbatim}

Moreover, Herzog and Hibi~\cite{HibiHerzog} proved that the projective dimension of $H_P$ is
the Dilworth number of $H_P$, i.e., the maximum length of an antichain of $H_P$.

\begin{verbatim}
i26 : pdim module HP == dilworthNumber P
o26 = true
\end{verbatim}

%


\end{document}